\documentclass{amsart}[11pt]
\usepackage{amssymb, mathrsfs, amsmath}
\newtheorem*{nota}{Remark}

\numberwithin{equation}{section}

\def\t{{\theta}}
\def\p{{\phi}}
\def\g{{\gamma}}

\def\la{{\langle}}
\def\k{{\kappa}}
\def\ra{{\rangle}}
\newtheorem{thm}{Theorem}[section]
\newtheorem{cor}[thm]{Corollary}

\newtheorem{prop}[thm]{Proposition}

\title{product formula for Jacobi polynomials, spherical harmonics and generalized Bessel function of dihedral type}

\begin{document}
\maketitle
\centerline{NIZAR DEMNI\footnote{SFB 701. Fakult\"at f\"ur mathematik, universit\"at Bielefeld. Germany. e-mail: demni@math.uni-bielefeld.de  
}} 
\begin{abstract}
We work out the expression of the generalized Bessel function of $B_2$-type derived in \cite{Dem1}. This is done using Dijskma and Koornwinder's product formula for Jacobi polynomials and the obtained expression is given by multiple integrals involving only a normalized modified Bessel function and two symmetric Beta distributions. 
We think of that expression as the major step toward the explicit expression of the Dunkl's intertwining operator $V_k$ in the $B_2$- invariant setting. Finally, we give in the same setting an explicit formula for the action of $V_k$ on a product of $|y|^{2\k}, \k \geq 0$ and the ordinary spherical harmonic $Y_{4m}(y) := |y|^{4m}\cos(4m\theta), y = |y|e^{i\theta}$. The obtained formula extends to all dihedral systems and it improves the one derived in \cite{X98}.  
\end{abstract}

\tableofcontents

\section{Motivation: Dunkl's intertwining operator}
The reader is referred to \cite{Dunkl} and references therein for an extensive background on the theory of Dunkl operators.     
The most challenging problem in this theory is an explicit expression for the action of the so-called Dunkl intertwining operator denoted in literature by $V_k$. The latter is a linear isomorphism of the space polynomials in several variables which is degree-preserving, conservative and which intertwines the commutative algebra of Dunkl operators and the one of partial derivatives. Few, yet relevant, works aiming to solve the above-mentioned problem were achieved and complicated formulae for the action of $V_k$ were supplied (\cite{Dunkl},
\cite{Dunkl0}, \cite{Xu}) except in the orthogonal case corresponding to products of copies of the group $\mathbb{Z}/2\mathbb{Z}$. In that case, $V_k$ is represented by means of a multivariate Beta distribution which reduces to a symmetric one (invariant under sign changes maps) in the W-invariant setting ($W$ is the reflections group), while for other nonorthogonal root systems multiple integrals with complicated integrands were obtained. Nevertheless, a relatively easy and elegant formula, compared to others, for the action of $V_k$ on $\langle x, \cdot \rangle$, where $x \in \mathbb{R}^2$ is fixed and $\langle\cdot,\cdot\rangle$ is the Euclidean scalar product, was given in \cite{X98} for the root system of type $B_2$. This root system is the most elementary example of non orthogonal root system and matches with the dihedral system of order $4$, denoted by ${\it I}_2(4)$ and corresponding to the symmetry group of the square. A noticeable coincidence is that while the proof supplied in \cite{X98} relies heavily on the product formula for Jacobi polynomials due to Dijskma and Koornwinder (\cite{Dij}), the expression derived in \cite{Dem1} for the generalized Bessel function of dihedral type, say $D_k^W$, of dihedral type also involves a product of Jacobi polynomials. Thus it is quite natural and interesting to use the product formula for Jacobi polynomials in order to seek a more elegant expression for $D_k^W$ in the $B_2$-type setting, aiming to investigate the connection between both works \cite{Dem1} and \cite{Dij} and to have a better insight into the action of $V_k$ on $B_2$-type invariant functions. 
In order to motivate the reader, we show below how the easy product formulae for cosine and sine functions allow, in some particular cases, to get an easy expression for $D_k^W$ of $B_2$-type. The general setting, corresponding to general Jacobi polynomials, is more subtle and one is far from getting (even from hoping) an easy expression for $D_k^W$ as in those particular cases. Nevertheless, the formula we supply involves only a modified Bessel function and two symmetric Beta distributions and simplifies considerably for some particular values of the multiplicity function, yielding as a by-product to the positivity of $D_k^W$.  
\subsection{Particular cases} 
Recall that even dihedral groups ${\it D}_2(2p), p \geq 2$ are the symmetry group of a $2p$-gone and that it contains two classes of reflections, with mirrors joining opposite vertices and mirrors joining the midpoints of opposite sides (\cite{Hum}). The following formula was derived in \cite{Dem1}: let $x= \rho e^{i\phi}, y = re^{i\theta}$ where $\rho,r \geq 0$ and $0 \leq \phi,\theta \leq \pi/(2p)$, then  

\begin{equation} \label{GBF}
D_k^W(x,y) = D_k^W(\rho,\phi,r,\theta) = c_{p,k}\left(\frac{2}{r\rho}\right)^{\gamma}
\sum_{j \geq 0}{\it I}_{2jp +\gamma}(\rho r)p_j^{l_1, l_0}(\cos(2p\phi))p_j^{l_1, l_0}(\cos(2p\theta))
\end{equation}
where $k=(k_0,k_1)$ is a positive-valued multiplicity function, $l_i = k_i - 1/2, i \in \{1,2\}$ are the index values, $\gamma = p(k_0+k_1)$ and ${\it I}_{\nu}, p_j^{l_1,l_0}$ are the modified Bessel function and the $j$-th (orthonormal) Jacobi polynomial respectively. A similar formula holds for odd dihedral systems $D_2(n), n\geq 1$ where one has to substitute in \eqref{GBF} $k_1=0, p=n, k_0 = k \geq 0$. The constant $c_{p,k}$ is such that $D_k^W(0,y) = 1$ for all $y$ and is easily computed as follows: recall the series expansion of the modified Bessel function:  

\begin{equation}\label{BessExp}
{\it I}_\nu (z) = \sum_{m=0}^\infty \frac{1}{m! \Gamma(m+\nu +1)} \left(\frac{z}{2}\right)^{2m+\nu}. 
\end{equation} 
Thus, for $\rho = 0$, the only non zero term in the double infinite series (after substitution of \eqref{BessExp}) corresponds to $j=m=0$, therefore 
\begin{equation*}
c_{p,k}^{-1} = \frac{1}{\Gamma(\gamma+1)}p_0^{l_1, l_0}(\cos(2p\phi))p_0^{l_1, l_0}(\cos(2p\theta)).
\end{equation*}
Since $p_0^{l_1,l_0}$ is constant and has unit-norm in $L^2([-1,1],(1-u)^{l_1}(1+u)^{l_0}du)$, then 
\begin{equation*}
\frac{1}{p_0^{l_1,l_0} (\cos(2p\theta))p_0^{l_1,l_0}(\cos(2p\phi))} = 2^{k_0+k_1}B(k_1+1/2,k_0+1/2)   
\end{equation*}
so that 
\begin{equation}\label{NormCons}
c_{p,k} = 2^{k_0+k_1}\frac{\Gamma(p(k_1+k_0)+1)\Gamma(k_1+1/2)\Gamma(k_0+1/2)}{\Gamma(k_0+k_1+1)}.
\end{equation}

The particular cases we mentioned above correspond to $k_0=k_1=0$ and $k_0=k_1=1$ and are referred to as the geometric cases. In the former, $P_j^{-1/2,-1/2}$ reduces to the orthonormal Tchebichef polynomial of the first kind defined by (\cite{Erd}) 
\begin{equation*}
T_j(\cos \theta) = \sqrt{\frac{2}{\pi}} \cos(j\theta),  j \neq 0, T_0(\cos \theta) = \sqrt{\frac{1}{\pi}}
\end{equation*}
so that $D_k^W$ simplifies to 
\begin{align}\label{D1}
D_0^W(\rho,\phi,r,\theta) &= \frac{c_p}{\pi}\left[I_0(\rho r) + 2 \sum_{j \geq 1} {\it I}_{2jp}(\rho r)\cos(2jp\phi)\cos(2jp\theta)\right] 
\\& =  {\it I}_0(\rho r) +  \sum_{j \in \mathbb{Z} \setminus \{0\}} {\it I}_{2jp}(\rho r)\cos(2jp\phi)\cos(2jp\theta) \nonumber
\\& =  \sum_{j \in \mathbb{Z}} {\it I}_{2jp}(\rho r)\cos(2jp\phi)\cos(2jp\theta)\nonumber.
\end{align}
since $c_p := c_{p,0} = \pi$ and  where the second equality follows from ${\it I}_{-m}(z) = {\it I}_m(z)$ for $m \in \mathbb{Z}$ (which is a consequence of $J_{-m}(z) = (-1)^mJ_m(z)$ and ${\it I}_{\nu}(z) = (-i)^{\nu}{\it J}_{\nu}(iz)$, see \cite{Wat}, p.15). Using the product formula $2\cos a \cos b = \cos(a+b) + \cos(a-b)$, one is led to compute (both series converge as we shall see)
\begin{align*}
\sum_{j \in \mathbb{Z}}{\it I}_{2jp}(\rho r)\cos(2jp(\phi \pm \theta)) &= \frac{1}{2}\sum_{j \in \mathbb{Z}}{\it I}_{2jp}(\rho r)\left[e^{i2jp(\phi \pm \theta)} + e^{-i2jp(\phi \pm \theta)}\right]
= \sum_{j \in \mathbb{Z}}{\it I}_{2jp}(\rho r)e^{i2jp(\phi \pm \theta)}. 
\end{align*}
Now, let $\xi = e^{i\pi/p}$ and $j \in \mathbb{Z}$, then it is an easy exercise to prove that 
\begin{equation}\label{Identity}
 \frac{1}{2p}\sum_{s=1}^{2p} \xi^{sj} = 
\left\{\begin{array}{lcr} 
1 & \textrm{if} & j \equiv 0[2p], \\
0 & \textrm{otherwise}, & 
\end{array}\right.
\end{equation}
for every $p \geq 2$. It follows that 
\begin{align*}
\sum_{j \in \mathbb{Z}}{\it I}_{2jp}(\rho r) e^{i2jp(\phi \pm \theta)} =  \frac{1}{2p}\sum_{s=1}^{2p}\sum_{j \in \mathbb{Z}}{\it I}_{j}(\rho r)e^{i(\phi \pm \theta)j}\xi^{js}.
\end{align*}
Now using the generating function for Bessel functions $($\cite{Wat} p.14$)$
\begin{equation}\label{Bessel}
e^{(t - 1/t)z/2} = \sum_{j \in \mathbb{Z}}{\it J}_j(z)t^j,\quad t \in \mathbb{C}\setminus \{0\},
\end{equation}
together with ${\it I}_{j}(z) = (-i)^{j}{\it J}_{j}(iz)$, one gets 
\begin{equation*}
\sum_{j \in \mathbb{Z}}{\it I}_{2jp}(\rho r)\cos(2jp(\phi \pm \theta)) = \frac{1}{2p}\sum_{s=1}^{2p}\sum_{j \in \mathbb{Z}}I_j(\rho r)e^{ij[(\phi \pm \theta)+ s\pi/p]} 
= \frac{1}{2p}\sum_{s=1}^{2p}e^{\rho r \cos[(\phi \pm \theta)+ s\pi/p]}.
\end{equation*}
As a result
\begin{equation*}
D_0^W(x,y) = \frac{1}{4p}\left[\sum_{s=1}^{2p}e^{\rho r \cos[(\phi + \theta)+ s\pi/p]} + \sum_{s=1}^{2p}e^{\rho r \cos[(\phi - \theta)+ s\pi/p]}\right]
\end{equation*}
which is in agreement with the definition  
\begin{equation*}
D_0^W(x,y) := \frac{1}{|W|}\sum_{w \in W}e^{\langle x,wy\rangle} = \frac{1}{4p}\sum_{w \in W}e^{\Re(x(\overline{wy}))} 
\end{equation*} 
and with the fact that the dihedral group contains $2p$ rotations $(x \mapsto x \xi^s, 1 \leq s \leq 2p)$ and $2p$ reflections $(x \mapsto \overline{x} \xi^s,, 1 \leq s \leq 2p)$. Note that for $p=2$, since $s\pi/2 \in \{\pi/2, \pi, 3\pi/2, 2\pi\}$, then easy computations transform $D_0^W$ to 
\begin{equation}\label{GBF0}
\frac{1}{4}\left[ \cosh(\rho r \cos(\phi + \theta)) + \cosh(\rho r \cos(\phi - \theta)) + \cosh(\rho r \sin(\phi + \theta)) + \cosh(\rho r \sin(\phi - \theta))\right].
\end{equation} 
Similarly, in the latter ($k_0=k_1=1$), $P_j^{1/2,1/2}$ reduces to the (orthonormal) Tchebichef polynomial of the second kind $U_j$ defined by (\cite{Erd})
\begin{equation*}
U_j(\cos(\theta)) = \sqrt{\frac{2}{\pi}}\frac{\sin((j+1)\theta)}{\sin \theta}, j \geq 0 .
\end{equation*}
Using the product formula $2\sin a\sin b = \cos(a-b) - \cos(a+b)$, one finds with $\omega(r,\theta):= r^{2p}\sin(2p\theta)$ that 
\begin{align*}
D_1^W(x,y) & = \frac{1}{\omega(r,\theta)\omega(\rho,\phi)}\sum_{s=1}^{2p}\sum_{j \in \mathbb{Z}}{\it I}_{j}(\rho r)[e^{i(\phi - \theta)j} - e^{i(\phi + \theta)j}]\xi^{js}
\\& =   \frac{1}{\omega(r,\theta)\omega(\rho,\phi)}\sum_{s = 1}^{2p}\left[e^{\rho r \cos[(\phi - \theta)+ s\pi/p)]} - e^{\rho r \cos[(\phi + \theta)+ s\pi/p)]}\right].
\end{align*} 
The last expression agrees with 
\begin{equation*}
D_1^W(x,y) =  \frac{1}{\omega(r,\theta)\omega(\rho,\phi)}\sum_{w \in W}\det(w)e^{\langle x,wy\rangle}
\end{equation*}
which follows from the shift principle (\cite{Dunkl}) or the reflection principle and the Doob's transform (\cite{Dem2}, \cite{Gra}). 

\subsection{Results} For general multiplicity values, the situation is far from being easy since for instance, $D_k^W$ of type $B$ is given by a multivariate series 
${}_0F_1$ defined via Jack polynomials (\cite{Dem2}) whose explicit expressions are not easy to write down (though note the recent progress in this direction) except in the few geometric cases. Nevertheless, we succeeded to express $D_k^W$ of $B_2$-type as a multiple integral involving only a normalized modified Bessel function and two symmetric Beta probability measures depending on $k_1,k_0$, that is we got rid of the dependence of $D_k^W$ on Jacobi polynomials. More explicitly, let  
 
\begin{equation*}
c_{2\phi,2\theta}(u,v) := \sqrt{\frac{1+u\cos 2\theta \cos 2\phi + v\sin2\theta\sin 2\phi}{2}}, \, (u,v) \in [-1,1]^2,
\end{equation*}  
and, for $t \in [0,1], 0 \leq z \leq 1/2$, let 
\begin{equation*}
K_{\gamma}(t,z) :=  \frac{1}{\sqrt{\pi}}\int_0^1\frac{q^{\gamma/2 -1}}{\sqrt{1-q}}  i_{\gamma/2}(t\sqrt{1-2zq}) dq,\, \gamma > 0,
\end{equation*}
$i_{\nu}$ being the normalized modified Bessel function\footnote{This function is slightly different from the spherical modified Bessel function which is valued $1$ at $0$.} defined by 
\begin{equation}\label{ModBesFun}
i_{\nu}(z) := \sum_{m=0}^{\infty}\frac{1}{\Gamma(m+\nu+1)m!} \left(\frac{z}{2}\right)^{2m}.
\end{equation} 
Denote also $\mu^{\nu}$ the symmetric Beta distribution
\begin{equation}
\mu^{\nu}(du) \propto (1-u^2)^{\nu-1/2}{\bf 1}_{[-1,1]}(u)du,\quad \nu > -1/2.
\end{equation}
Then, our first result may be stated as follows:

\begin{prop}\label{Prop1}
\begin{align}\label{ObF}
&D_k^W(x,y) = \Gamma\left(\frac{\gamma+1}{2}\right) \int\int \frac{1+\cos^2(\gamma \pi/2)}{2} i_{(\gamma-1)/2}\left(\rho r c_{2\phi,2\theta}(u,v)\right) +\sin^2(\gamma \pi/2) \nonumber
\\&\int_0^{1} {\it I}_0(\rho rt) \partial_t\left\{t^{\gamma}\left[\gamma K_{\gamma}\left(\rho rt, \frac{c_{2\phi,2\theta}^2(u,v)}{2}\right) 
+ \frac{\gamma+1}{\gamma} \partial_zK_{\gamma} \left(\rho rt, \frac{c_{2\phi,2\theta}^2(u,v)}{2}\right)\right]\right\} dt.
\mu^{l_1}(du)\mu^{l_0}(dv),
\end{align}
The extension to $\gamma=0$ is performed via weak limit and one recovers \eqref{GBF0}. 
\end{prop}

\begin{nota}
\eqref{ObF} is relatively easy and explicit regarding the existing formulae (\cite{Dunkl0}, \cite{Xu}). Besides, we think of it as a major step toward an explicit expression for the action of $V_k$ on $B_2$-invariant functions since Bessel functions enjoy a huge number of nice properties. However we do not know how to come to $V_k$ due mainly to the last term in the integrand. 
\end{nota}

\begin{cor}
For even values of $\gamma \neq 0$,
\begin{align*}
D_k^W(x,y) = \Gamma\left(\frac{\gamma+1}{2}\right) \int\int i_{(\gamma-1)/2}\left(\rho r c_{2\phi,2\theta}(u,v)\right)\mu^{l_1}(du)\mu^{l_0}(dv) 
\end{align*}
which is obviously positive. 
\end{cor}

Our second result supplies for all dihedral systems, independently from the first one, an explicit formula for the action of $V_k$ on a product of $|y|^{2\k}, \k \geq 0$ and the ordinary spherical harmonic $Y_{m}(y) = |y|^{m}\cos(m\theta), y = |y|e^{i\theta}$. For sake of simplicity and coherency with the first result, we only write down the formula for the $D_2(4)$-type root system:  
\begin{prop}\label{Prop2}
For positive integers $\k,m$ and $y = |y|e^{i\theta}$, one has 
\begin{equation*}
V_k \left[|\cdot |^{2\k} Y_{4m} \right](y)  =   |y|^{2\k+4m} \sum_{0 \le 2j \le \k}  \frac{ b_{m,j+m}\k!\Gamma(4m+\k+1) } {(\k-2j)! \Gamma(4m+2j +\k + \g +1)}p_{j+m}^{l_1, l_0}(\cos(4\t)),
\end{equation*}
where $b_{m,j+m}$ is the $m$-th coefficient in the expansion of $P_{j+m}^{l_1, l_0}(\cos(4\t))$ as a finite linear combination of $\cos(4m\t)$. 
\end{prop}

\begin{nota}
Proposition \eqref{Prop2} shows the action of $V_\k$ on polynomials that are invariant under the group $B_2$. Such results have been studied and can be derived after some manipulations from \cite{X98}, but the formula obtained here is more explicit.  We didn't attempt to compute the special cases in which the explicit formula of $b_{j,m}$ are known. This is of interest only if the formula can be further simplified or the sum in the proposition can be summed up.  Note also that our formula resembles the one displayed in Theorem.2 p.13 in \cite{Dunkl2}.
\end{nota}


\section{$B_2$-type root system: proof of Proposition \eqref{Prop1}}
As indicated above, this section consists of six subsections, the title of each one indicates its content, aiming to orient the reader through the overwhelming yet tricky computations. Before coming into them, we want to inform the reader that usual operations performed on integrals (exchange of integration's order or limiting operations under integral signs) are easily justified due to the compactness of integration domains and to the smoothness of the involved functions.     

\subsection{Product formula: a first transformation}
The generalized Bessel function displayed in \eqref{GBF} reads in the $B_2$-type setting
\begin{align*}
D_k^W(x,y) &=  c_{2,k}\left(\frac{2}{r\rho}\right)^{\gamma}\sum_{j \equiv 0[4]}{\it I}_{j +\gamma}(\rho r)p_{j/4}^{l_1, l_0}(\cos(4\phi))p_{j/4}^{l_1, l_0}(\cos(4\theta))
\end{align*}
where 
\begin{equation*}
c_{2,k} = \frac{2^{3\gamma/2}}{\sqrt{\pi}}\Gamma\left(\frac{\gamma+1}{2}\right)\Gamma(k_1+1/2)\Gamma(k_0+1/2)
\end{equation*}
by Gauss duplication's formula. Now, recall the product formula for Jacobi polynomials (\cite{Dij}\footnote{We write the formula using orthonormal Jacobi polynomials and symmetric Beta probability measures for later purposes, yet the Gegenbauer polynomial is not normalized.}, \cite{Xu} p.424\footnote{There is an erratum in the constant term in front of the integral.})
\begin{equation*}
c_{\alpha,\beta} p_j^{\alpha,\beta}(\cos2\phi)p_j^{\alpha,\beta}(\cos2\theta) = (2j+\alpha+\beta+1)\int\int C_{2j}^{\alpha+\beta+1}(z_{\theta}(u,v))\mu^{\alpha}(du)\mu^{\beta}(dv)
\end{equation*}  
where $\Re(\alpha),\Re(\beta) > -1/2, \phi,\theta \in [-\pi/2,\pi/2]$, 
\begin{equation*}
c_{\alpha,\beta} = 2^{\alpha+\beta+1}\frac{\Gamma(\alpha+1)\Gamma(\beta+1)}{\Gamma(\alpha+\beta+1)},
\end{equation*}
\begin{equation}\label{Term}
z_{\phi,\theta}(u,v) = u\cos \theta\cos\phi + v \sin\theta\sin\phi
\end{equation}
and $\mu^{\alpha}$ is the symmetric Beta probability measure whose density is given by
\begin{equation*}
\mu^{\alpha}(du) = \frac{\Gamma(\alpha+1)}{\sqrt{\pi}\Gamma(\alpha+1/2)}(1-u^2)^{\alpha-1/2} {\bf 1}_{[-1,1]}(u)du, \quad \alpha > -1/2.
\end{equation*}
Specializing the product formula to $\alpha = l_1 = k_1-1/2 > -1/2,\beta = l_0 = k_0-1/2 > -1/2$ gives 
\begin{align*}
D_k^W(x,y) &=  \frac{2^{\gamma}}{\sqrt{\pi}}\Gamma\left(\frac{\gamma+1}{2}\right)\Gamma\left(\frac{\gamma}{2}\right)\left(\frac{2}{r\rho}\right)^{\gamma}
\sum_{j \equiv 0[4]}{\it I}_{j +\gamma}(\rho r)(j+\gamma)/2\int\int C_{j/2}^{k_1+k_0}(z_{2\phi,2\theta}(u,v))\mu^{l_1}(du)\mu^{l_0}(dv)
\\& =  \Gamma(\gamma)\left(\frac{2}{r\rho}\right)^{\gamma}
 \sum_{j \equiv 0[4]}{\it I}_{j +\gamma}(\rho r)(j+\gamma)\int\int C_{j/2}^{k_1+k_0}(z_{2\phi,2\theta}(u,v))\mu^{l_1}(du)\mu^{l_0}(dv)
\end{align*}
by Gauss's duplication formula. Now, a useful by-product of the product formula is the following (Theorem 2.2. in \cite{Dij} specialized to $\alpha= \beta = \lambda - 1/2, \lambda > 0$)
\begin{equation*}
C_j^{\lambda}(z) = \int C_{2j}^{2\lambda}\left(\sqrt{\frac{1+z_{2\phi,2\theta}}{2}}w\right) \mu^{\lambda-1/2}(dw)
\end{equation*}
and yields with $\lambda  = \gamma/2$
\begin{align*}
&D_k^W(x,y) = \Gamma(\gamma)\left(\frac{2}{r\rho}\right)^{\gamma}\sum_{j \equiv 0[4]}{\it I}_{j +\gamma}(\rho r)(j+ \gamma)\int \int\int C_{j}^{\gamma}\left(Z\right) 
\mu^{l_1}(du)\mu^{l_0}(dv)\mu^{(\gamma-1)/2}(dw)
\end{align*}
where we set
\begin{equation*}
Z = Z_{2\phi,2\theta}(u,v,w) := \sqrt{\frac{1+z_{2\phi,2\theta}(u,v)}{2}}w.
\end{equation*}
Finally, using (\ref{Identity}), one gets: 
\begin{align}\label{Dunkl}
&D_k^W(x,y) = \frac{\Gamma(\gamma)}{4} \left(\frac{2}{r\rho}\right)^{\gamma}\sum_{s=1}^{4} \sum_{j \geq 0}{\it I}_{j +\gamma}(\rho r)(j+ \gamma)\xi^{js}
\int \int\int C_{j}^{\gamma}\left(Z\right) \mu^{l_1}(du)\mu^{l_0}(dv)\mu^{(\gamma-1)/2}(dw).
\end{align}

\subsection{An auxiliary formula}
Since 
\begin{equation}\label{BF}
|C_j^{\gamma}(Z)| \leq |C_j^{\gamma}(1)| = \frac{(2\gamma)_j}{j!}, 
\end{equation}
and since
\begin{equation*}
\sum_{j \geq 0}\frac{(2\gamma)_j}{j!} {\it I}_{j+\gamma}(\rho r) \, < \, \infty,
\end{equation*} 
then Fubini's Theorem applies and one is led to compute 
\begin{equation*}
\sum_{j \geq 0}(j+ \gamma){\it I}_{j +\gamma}(\rho r) C_{j}^{\gamma}\left(Z\right)\xi^{js}
\end{equation*}
for fixed $(u,v,w) \in ]-1,1[^3$. Note that for $s= 2,4$, one has 
\begin{equation}\label{Form}
 \sum_{j \geq 0}(\pm 1)^j(j+ \gamma){\it I}_{j +\gamma}(\rho r)C_{j}^{\gamma}(Z) = \frac{(\rho r)^{\gamma}}{2^{\gamma}\Gamma(\gamma)}e^{\pm \rho r Z}
 \end{equation}
by formula 5.13.3.3. p.712 in \cite{Bry}, yet we did not find any similar result to 
\begin{equation*}
\sum_{j \geq 0}(\pm i)^j(j+ \gamma){\it I}_{j +\gamma}(\rho r)C_{j}^{\gamma}(Z).
\end{equation*}
However, note that for strictly positive integer values of $\gamma$, \eqref{Form} may be written as 
\begin{equation*}
 \sum_{j \geq 0}(\pm 1)^j(j+ \gamma){\it I}_{j +\gamma}(\rho r)C_{j}^{\gamma}(Z) = \frac{(\pm 1)^{\gamma}}{2^{\gamma}\Gamma(\gamma)}D_Z^{\gamma}\left[e^{\pm \rho r Z}\right]. 
 \end{equation*}
In fact, we can derive a more general similar result for  
\begin{equation*}
\sum_{j \geq 0}(j+ \gamma){\it I}_{j +\gamma}(\rho r)C_{j}^{\gamma}(Z)e^{ijq}
\end{equation*}
for any real number $q$ provided that $\gamma \geq 1$ is an integer (then extend it to strictly positive values). In fact, from p.32 and p.44 in \cite{Ober}, one gets
\begin{align*}
2\sum_{j \geq 0}{\it I}_j(z)e^{ij\zeta} &=  2\sum_{j \geq 0}{\it I}_j(z)[\cos(j\zeta)+i \sin(j\zeta)]
\\& = \exp(z\cos \zeta) + {\it I}_0(\zeta) + i \sin \zeta \int_0^z\exp(-t\cos \zeta){\it  I}_0(t)dt
\end{align*}
for all real numbers $z,\zeta$ and from 11.1.2 (18) p.235 in \cite{Erd}, the following holds
\begin{equation*}
(j+\gamma)C_j^{\gamma}(Z) = \frac{1}{2^{\gamma-1}\Gamma(\gamma)} D_Z^{\gamma}(T_{j+\gamma}(Z)) = \frac{1}{2^{\gamma-1}\Gamma(\gamma)}
D_Z^{\gamma}[\cos((j+\gamma)\arccos Z)], 
\end{equation*}
where $D_Z$ stands for the derivative operator. Moreover by the bound \eqref{BF}, one has
 \begin{align*}
\sum_{j \geq 0}(j+ \gamma){\it I}_{j +\gamma}(\rho r) C_{j}^{\gamma}\left(Z\right)e^{ijq} &= \frac{1}{2^{\gamma-1}\Gamma(\gamma)}
\sum_{j \geq 0}{\it I}_{j +\gamma}(\rho r) D_Z^{\gamma}[\cos((j+\gamma)\arccos Z)]e^{ijq}
\\& = \frac{\xi^{-s\gamma}}{2^{\gamma-1}\Gamma(\gamma)} D_Z^{\gamma}\sum_{j \geq \gamma}{\it I}_{j}(\rho r) \cos(j\arccos Z)e^{ijq},
\end{align*}
and since $\cos(j\arccos Z)$ is a (Tchbeycheff) polynomial of degree $j$, then
\begin{align*}
\sum_{j \geq 0}(j+ \gamma){\it I}_{j +\gamma}(\rho r) C_{j}^{\gamma}\left(Z\right)e^{ijq}&=
\frac{\xi^{-s\gamma}}{2^{\gamma-1}\Gamma(\gamma)} D_Z^{\gamma}\sum_{j \geq 0}{\it I}_{j}(\rho r) \cos(j\arccos Z)e^{ijq}
\\& = \frac{\xi^{-s\gamma}}{2^{\gamma}\Gamma(\gamma)}D_Z^{\gamma}\sum_{j \geq 0}{\it I}_{j}(\rho r) [e^{ija} + e^{-ija}]e^{ijq},
\end{align*}
where we set $a := \cos Z$. 

\subsection{Odd and even values of $\gamma$} Easy computations using the last above equality yield
 \begin{align*}
\sum_{j \geq 0}(\pm i)^j(j+ \gamma){\it I}_{j +\gamma}(\rho r)C_{j}^{\gamma}(Z) &=  \frac{(\mp i)^{\gamma}}{2^{\gamma}\Gamma(\gamma)}D_Z^{\gamma} 
\left[\cosh(\rho r\sqrt{1-Z^2})] \pm i Z \int_0^{\rho r}\cosh(t\sqrt{1-Z^2})I_0(t)dt \right]
\end{align*}
for $s=1,3$. If $\gamma$ is even, then $i^{\gamma} = (-i)^{\gamma}$ so that
\begin{align}\label{Aux1}
(+)\sum_{j \geq 0}(j+ \gamma){\it I}_{j +\gamma}(\rho r) C_{j}^{\gamma}\left(Z\right)(\pm i)^{js} & =  \frac{(-1)^{\gamma/2}}{2^{\gamma-1}\Gamma(\gamma)}D_Z^{\gamma} 
[\cosh(\rho r\sqrt{1-Z^2})],
\end{align}
while if $\gamma $ is odd then $(-i)^{\gamma} = -i (-i)^{\gamma-1} =  -i (i)^{\gamma-1} = - (i)^{\gamma}$ therefore
\begin{align}\label{Aux2}
(+)\sum_{j \geq 0}(j+ \gamma){\it I}_{j +\gamma}(\rho r) C_{j}^{\gamma}\left(Z\right)(\pm i)^{js} =  \frac{(-1)^{(\gamma-1)/2}}{2^{\gamma-1}\Gamma(\gamma)}
D_Z^{\gamma}[Z \int_0^{\rho r}\cosh(t\sqrt{1-Z^2})I_0(t)dt].
\end{align}
Above we used the symbol $(+)$ to indicate that we sum both series. For $s=2,4$, one recovers \eqref{Form} and both series displayed there contribute to 
\begin{equation*}
(+) \sum_{j \geq 0}(\pm 1)^j(j+ \gamma){\it I}_{j +\gamma}(\rho r)C_{j}^{\gamma}(Z) = \frac{(\rho r)^{\gamma}}{2^{\gamma-1}\Gamma(\gamma)}\cosh (\rho r Z).
 \end{equation*}

Let $i_{\nu}$ denote the normalized modified Bessel function already defined in \eqref{ModBesFun}, it has the Poisson integral representation for $\nu > -1/2$ (\cite{Erd} p.81):
\begin{equation*}
i_{\nu}(z) = \frac{1}{\sqrt{\pi}\Gamma(\nu+1/2)}\int_{-1}^1e^{zw}(1-w^2)^{\nu-1/2}dw = \frac{1}{\Gamma(\nu+1)}\int_{-1}^1\cosh(zw)\mu^{\nu}(dw).
\end{equation*} 
Therefore, one gets with $\nu = (\gamma-1)/2$, 
\begin{align}\label{Int1}
&\frac{2^{\gamma}\Gamma(\gamma)}{(\rho r)^{\gamma}}\int_{-1}^1 (+) \sum_{j \geq 0}(\pm 1)^j(j+ \gamma){\it I}_{j +\gamma}(\rho r)C_{j}^{\gamma}(Z)\mu^{(\gamma-1)/2}(dw)= 
2\Gamma\left(\frac{\gamma+1}{2}\right) i_{(\gamma-1)/2}[\rho rc_{2\phi,2\theta}(u,v)],
\end{align} 
 where we set
\begin{equation*}
 Z = \sqrt{\frac{1+z_{2\phi,2\theta}(u,v)}{2}}w := c_{2\phi,2\theta}(u,v) w.
\end{equation*}
With regard to \eqref{Aux1}, one needs then to integrate 
\begin{align}\label{Int2}
\int_{-1}^1 D_Z^{\gamma} [\cosh(\rho r\sqrt{1-Z^2})]_{|Z = c_{2\phi,2\theta}(u,v) w} \mu^{(\gamma-1)/2}(dw), 
 \end{align}
for even $\gamma \geq 2$.

\subsection{A positive-definite function} 
One tricky way to compute \eqref{Int2} is to use the positive-definiteness of $Z \mapsto \cosh(t \sqrt{1-Z^2})$ for all real $t$. 
In fact, the following Bochner representation holds and is easily derived from 6.645.3 in \cite{Grad} by an analytic continuation
\begin{align*}
\cosh(t\sqrt{1-Z^2}) &= \cos (tZ) + \frac{t}{2}\int_{-1}^1 e^{itZq}\frac{{\it I}_1(t\sqrt{1-q^2})}{\sqrt{1-q^2}}dq :=\int e^{itZq}\nu_t(dq),
\end{align*}
where $\nu_t$ is the symmetric measure
\begin{equation*}
\nu_t(dq) = \frac{1}{2}\left[\delta_1(dq)+\delta_{-1}(dq)\right] + \frac{t}{2}\frac{{\it I}_1(t\sqrt{1-q^2})}{\sqrt{1-q^2}}{\bf 1}_{]-1,1[}(q)dq.
\end{equation*}
Thus, one has for even $\gamma$ 
\begin{align*}
D_Z^{\gamma} \cosh(\rho r\sqrt{1-Z^2}) = (-1)^{\gamma/2}(\rho r)^{\gamma}  \int q^{\gamma}e^{itZq}\nu_{\rho r}(dq).
\end{align*} 
Using Fubini's Theorem, it follows that 
\begin{align*}
& (-1)^{\gamma/2} \int_{-1}^{1} D_Z^{\gamma}[\cosh(\rho r\sqrt{1-Z^2})]_{|Z = c_{2\phi,2\theta}(u,v) w} \mu^{(\gamma-1)/2}(dw) = 
\\& (\rho r)^{\gamma}\int q^{\gamma} \left(\int_{-1}^{1}e^{i\rho r c_{2\phi,2\theta}(u,v)wq}\mu^{(\gamma-1)/2}(dw)\right)\nu_{\rho r}(dq)
\end{align*}
and  the integral between brackets is nothing but the spherical Bessel function $j_{(\gamma-1)/2}$ defined by $j_{\nu}(z) :=  i_{\nu}(iz)$ so that:
\begin{align*}
 (-1)^{\gamma/2} \int_{-1}^{1} D_Z^{\gamma}[\cosh(\rho r\sqrt{1-Z^2})]_{|Z = c_{2\phi,2\theta}(u,v) w} &\mu^{(\gamma-1)/2}(dw) =(\rho r)^{\gamma} \Gamma\left(\frac{\gamma+1}{2}\right) \\&
 \int q^{\gamma}  j_{(\gamma-1)/2}[c_{2\phi,2\theta}(u,v)\rho rq] \nu_{\rho r}(dq).
\end{align*}
With regard to \eqref{Aux1} and \eqref{GBF0}, 
\begin{align}\label{Int21}
\frac{2^{\gamma}\Gamma(\gamma)}{(\rho r)^{\gamma}}\int_{-1}^1(+)\sum_{j \geq 0}(j+ \gamma){\it I}_{j +\gamma}(\rho r) C_{j}^{\gamma}\left(c_{2\phi,2\theta}(u,v) w\right)&(\pm i)^{js}\mu^{(\gamma-1)/2}(dw) =  \nonumber
2\Gamma\left(\frac{\gamma+1}{2}\right) \\& \int q^{\gamma}  j_{(\gamma-1)/2}(c_{2\phi,2\theta}(u,v)\rho rq)\nu_{\rho r}(dq).
\end{align}
Now, using the expansions of $j_{\gamma-1/2}$ and ${\it I}_1$, one has 
\begin{align*}
&\frac{\rho r}{2}\int_{-1}^1q^{\gamma}  j_{(\gamma-1)/2}(c_{2\phi,2\theta}(u,v)\rho r q)\frac{{\it I}_1(\rho r\sqrt{1-q^2})}{\sqrt{1-q^2}} dq = 
\sum_{j,m\geq 0}\frac{(-1)^j}{\Gamma(j+(\gamma+1)/2)j!} \frac{1}{(m+1)!m!} 
\\& \left(\frac{c_{2\phi,2\theta}(u,v) \rho r}{2}\right)^{2j} \left(\frac{\rho r}{2}\right)^{2m+2} \int_0^1q^{(\gamma-1)/2+j}(1-q)^{m}dq
\\& = \sum_{j,m\geq 0}\frac{(-1)^j}{\Gamma(j+(\gamma-1)/2 + m+2)j!} \frac{1}{(m+1)!}\left(\frac{c_{2\phi,2\theta}(u,v) \rho r}{2}\right)^{2j} \left(\frac{\rho r}{2}\right)^{2m+2}
\\& = \sum_{j \geq 0}\frac{(-1)^j}{j!}\left(\frac{c_{2\phi,2\theta}(u,v) \rho r}{2}\right)^{2j}\sum_{m \geq 1}\frac{1}{\Gamma(j+(\gamma-1)/2 + m+1)m!}\left(\frac{\rho r}{2}\right)^{2m}
\\& = \sum_{j \geq 0}\frac{(-1)^j}{j!}\left(\frac{c_{2\phi,2\theta}(u,v) \rho r}{2}\right)^{2j}\left[\left(\frac{2}{\rho r}\right)^{(\gamma-1)/2 + j}{\it I}_{j+(\gamma-1)/2}(\rho r) - 
\frac{1}{\Gamma((\gamma-1)/2+j+1)}\right]
\\& = \left(\frac{2}{\rho r}\right)^{(\gamma-1)/2}\sum_{j \geq 0}\frac{(-1)^j}{j!}\left(\frac{[c_{2\phi,2\theta}(u,v)]^2\rho r}{2}\right)^{j} {\it I}_{j+(\gamma-1)/2}(\rho r) - j_{(\gamma-1)/2}(c_{2\phi,2\theta}(u,v)\rho r).
\end{align*}
Using the formula (11.4. p.694 in \cite{Bry})
\begin{equation*}
\sum_{m \geq 0}\frac{z^k}{k!} {\it I}_{a+k}(x) = \left(1+\frac{2z}{x}\right)^{-a/2}{\it I}_a\left(x\sqrt{1+\frac{2z}{x}}\right),\, |2z| < |x|,
\end{equation*}
one finally sees that \eqref{Int21} simplifies to ($|c_{2\phi,2\theta}(u,v)| < 1$ for almost all $(u,v)$)
\begin{align*} 
\left(\frac{2}{\rho r\sqrt{1-c^2_{2\theta}(u,v)}}\right)^{(\gamma-1)/2}{\it I}_{(\gamma-1)/2}(\rho r\sqrt{1-c^2_{2\theta}(u,v)}) = i_{(\gamma-1)/2}(\rho r\sqrt{1-c^2_{2\theta}(u,v)}).
\end{align*}

Finally 
\begin{align}\label{Int2111}
\frac{2^{\gamma}\Gamma(\gamma)}{(\rho r)^{\gamma}}\int_{-1}^1(+)\sum_{j \geq 0}(j+ \gamma)&{\it I}_{j +\gamma}(\rho r) C_{j}^{\gamma}\left(c_{2\phi,2\theta}(u,v) w\right)(\pm i)^{js}\mu^{(\gamma-1)/2}(dw) =  \nonumber
2\Gamma\left(\frac{\gamma+1}{2}\right)\\& \cos^2\left(\frac{\gamma \pi}{2}\right)i_{(\gamma-1)/2}(\rho r\sqrt{1-c^2_{2\theta}(u,v)})
\end{align}
which extends to all strictly positive integer values of $\gamma$ as being zero at odd values. Now, let $\gamma$ be an odd positive integer, then by the virtue of \eqref{Aux2}, it remains only to integrate
\begin{align}\label{Int3}
\int_{-1}^1 D_Z^{\gamma} [Z\cosh(t\sqrt{1-Z^2})]_{|Z = c_{2\phi,2\theta}(u,v)w}\mu^{(\gamma-1)/2}(dw) 
 \end{align}
for fixed $u,v \in [-1,1]$ and $t \in [0,\rho r]$. Using Leibniz's rule, one gets:
\begin{equation*}
D_Z^{\gamma}[Z \cosh(t\sqrt{1-Z^2})] = Z D_Z^{\gamma}[\cosh(t\sqrt{1-Z^2})] + \gamma D_Z^{\gamma-1}[\cosh(t\sqrt{1-Z^2})],
\end{equation*}
so that an integration by parts transforms the second integral in \eqref{Int3} to
\begin{align*}
&\int_{-1}^1 D_Z^{\gamma}[Z\cosh(t\sqrt{1-Z^2})]_{|Z = c_{2\phi,2\theta}(u,v)w}\mu^{(\gamma-1)/2}(dw) = 
\\& \frac{1}{\gamma}\int_{-1}^1  D_Z^{\gamma+1}[\cosh(t\sqrt{1-Z^2})]_{|Z = c_{2\phi,2\theta}(u,v)w}\mu^{(\gamma+1)/2}(dw)
\\& + \gamma \int_{-1}^1 D_Z^{\gamma-1}[\cosh(t\sqrt{1-Z^2})]_{|Z = c_{2\phi,2\theta}(u,v)w}\mu^{(\gamma-1)/2}(dw).
\end{align*}
Hence, \eqref{Int3} is almost similar to \eqref{Int2} so that one gets: 
\begin{align*}
(-1)^{(\gamma-1)/2}\int_{-1}^1 D_Z^{\gamma+1}&[\cosh(t\sqrt{1-Z^2})]_{|Z = c_{2\phi,2\theta}(u,v)w}\mu^{(\gamma+1)/2}(dw) 
=-\Gamma\left(\frac{\gamma+3}{2}\right) t^{\gamma+1}\\& \int q^{\gamma+1}  j_{(\gamma+1)/2}(c_{2\phi,2\theta}(u,v)tq)\nu_t(dq)
\end{align*}
 and 
\begin{align*}
(-1)^{(\gamma-1)/2}\int_{-1}^{1} D_Z^{\gamma-1}&[\cosh(t\sqrt{1-Z^2})]_{|Z = c_{2\phi,2\theta}(u,v)w} \mu^{(\gamma-1)/2}dw =  \Gamma\left(\frac{\gamma+1}{2}\right)t^{\gamma-1}
\\& \int q^{\gamma-1}  j_{(\gamma-1)/2}(c_{2\phi,2\theta}(u,v)tq)\nu_t(dq).
\end{align*}

Also, similar computations as done before show that
\begin{align*}
\int q^{\gamma+1}  j_{(\gamma+1)/2}(c_{2\phi,2\theta}(u,v)tq)\nu_t(dq) &= \left(\frac{2}{t}\right)^{\gamma/2} 
\sum_{j \geq 0}\frac{(-1)^j\Gamma((\gamma/2)+j+1)}{\Gamma((\gamma+1)/2+j+1)j!}\left(\frac{t}{2}c_{2\phi,2\theta}^2(u,v)\right)^j {\it I}_{\gamma/2+j}(t)
\end{align*}
and using the derivation rule $[z^{\nu+1}{\it I}_{\nu+1}(z)]' = z^{\nu+1}{\it I}_{\nu}(z)$ (\cite{Wat}) with $\nu = \gamma/2 + j$, it then follows that
\begin{align}\label{Sum1}
&-\int (tq)^{\gamma+1}  j_{(\gamma+1)/2}(c_{2\phi,2\theta}(u,v)tq)\nu_t(dq) =  -2^{\gamma/2}\frac{d}{dt}
 \sum_{j \geq 0}\frac{(-1)^j\Gamma((\gamma/2)+j+1)}{\Gamma((\gamma+1)/2+j+1)j!}\left(\frac{c_{2\phi,2\theta}^2(u,v)}{2}\right)^{j} \nonumber
 \\& t^{\gamma/2+j+1} {\it I}_{\gamma/2+j+1}(t) =  2^{\gamma/2} \frac{d}{dt}\sum_{j \geq 1}\frac{(-1)^j\Gamma((\gamma/2)+j)}{\Gamma((\gamma+1)/2+j)(j-1)!}\left(\frac{c_{2\phi,2\theta}^2(u,v)}{2}\right)^{j-1}t^{\gamma/2+j} {\it I}_{\gamma/2+j}(t) \nonumber 
\\&=  2^{\gamma/2} \frac{d}{dt}\frac{d}{dz}\sum_{j \geq 0}\frac{(-1)^j\Gamma((\gamma/2)+j)}{\Gamma((\gamma+1)/2+j)j!}z^{j} t^{\gamma/2+j} {\it I}_{\gamma/2+j}(t)_{|z= c_{2\phi,2\theta}^2(u,v)/2}.
\end{align}
Substituting $\gamma+1$ by $\gamma-1$, one gets 
\begin{align*}
\int q^{\gamma-1}  j_{(\gamma-1)/2}(c_{2\phi,2\theta}(u,v)tq)\nu_t(dq) &= \left(\frac{2}{t}\right)^{\gamma/2-1} 
\sum_{j \geq 0}\frac{(-1)^j\Gamma(\gamma/2+j)}{\Gamma((\gamma+1)/2+j)j!}\left(\frac{t}{2}c_{2\phi,2\theta}^2(u,v)\right)^j {\it I}_{\gamma/2+j-1}(t)
\end{align*}
so that
\begin{align}\label{Sum2}
\int (tq)^{\gamma-1}  j_{(\gamma-1)/2}(c_{2\phi,2\theta}(u,v)tq)\nu_t(dq) &= 2^{\gamma/2-1} \frac{d}{dt}
\sum_{j \geq 0}\frac{(-1)^j\Gamma(\gamma/2+j)}{\Gamma((\gamma+1)/2+j)j!}\left(\frac{c_{2\phi,2\theta}^2(u,v)}{2}\right)^j \nonumber \\& t^{\gamma/2+j}{\it I}_{\gamma/2+j}(t).
\end{align}
Set 
\begin{align*}
A_{\gamma}(t,z) & := \sum_{j \geq 0}\frac{(-1)^j\Gamma((\gamma/2)+j)}{\Gamma((\gamma+1)/2+j)j!}z^{j} t^{\gamma/2+j} {\it I}_{\gamma/2+j}(t)
\\& =  \frac{t^{\gamma}}{2^{\gamma/2}}\sum_{j \geq 0}\frac{(-1)^j\Gamma((\gamma/2)+j)}{\Gamma((\gamma+1)/2+j)j!}\left(\frac{zt^2}{2}\right)^{j} i_{\gamma/2+j}(t)
:= \frac{t^{\gamma}}{2^{\gamma/2}}K_{\gamma}(t,z),
\end{align*}
for $\gamma > 0$ and $0 < z < 1/2$, then a glance at \eqref{Sum1} and \eqref{Sum2} shows that 

\begin{align*}
(-1)^{(\gamma-1)/2}\int_{-1}^1 &D_Z^{\gamma} [Z\cosh(t{1-Z^2})]_{|Z = c_{2\phi,2\theta}(u,v)w}\mu^{(\gamma-1)/2}(dw) = \frac{1}{2} \Gamma\left(\frac{\gamma+1}{2}\right)
\\&\partial_t\left\{t^{\gamma}\left[\gamma K_{\gamma}\left(t, \frac{c_{2\phi,2\theta}^2(u,v)}{2}\right) + \frac{\gamma+1}{\gamma} \partial_zK_{\gamma} \left(t, \frac{c_{2\phi,2\theta}^2(u,v)}{2}\right)\right]\right\}.
\end{align*}

Therefore,  with regard to \eqref{Aux2}, one has for odd $\gamma \geq 1$
\begin{align}\label{Int31}
&\frac{2^{\gamma}\Gamma(\gamma)}{(\rho r)^{\gamma}} \int_{-1}^1 \sum_{s=1}^4\sum_{j \geq 0}(j+ \gamma){\it I}_{j +\gamma}(\rho r) C_{j}^{\gamma}\left(c_{2\phi,2\theta}(u,v) w\right)\xi^{js}\mu^{(\gamma-1)/2}(dw) = \Gamma\left(\frac{\gamma+1}{2}\right)\nonumber 
\\& \sin^2\left(\frac{\gamma \pi}{2}\right) \int_0^{1} {\it I}_0(\rho rt) \partial_t\left\{t^{\gamma}\left[\gamma K_{\gamma}\left(\rho rt, \frac{c_{2\phi,2\theta}^2(u,v)}{2}\right) + \frac{\gamma+1}{\gamma} \partial_zK_{\gamma} \left(\rho rt, \frac{c_{2\phi,2\theta}^2(u,v)}{2}\right)\right]\right\} dt,
\end{align}
which extends to all strictly positive integer values of $\gamma$ as the zero function at even values. Unfortunately, we did not succeed to derive any more easier expression for $A_{\gamma}(t,z)$ (may be this is not possible), nevertheless noting that 
\begin{equation*}
\frac{\sqrt{\pi}\Gamma((\gamma/2)+j)}{\Gamma((\gamma+1)/2+j)} = \int_0^1 \frac{q^{\gamma/2 + j-1}}{\sqrt{1-q}} dq 
\end{equation*}   
for $\gamma > 0$, then
\begin{align*}
A_{\gamma}(t,z) &= \frac{t^{\gamma/2}}{\sqrt{\pi}} \int_0^1\frac{q^{\gamma/2 -1}}{\sqrt{1-q}}\sum_{j \geq 0}\frac{(-1)^j\Gamma((\gamma/2)+j)}{\Gamma((\gamma+1)/2+j)j!}(qtz)^{j}  {\it I}_{\gamma/2+j}(t)dq  
\\& = \frac{t^{\gamma/2}}{\sqrt{\pi}} \int_0^1\frac{q^{\gamma/2 -1}}{\sqrt{1-q}}\frac{1}{\sqrt{1-2zq}^{\gamma/2}} {\it I}_{\gamma/2}(t\sqrt{1-2zq})dq
\\& = \frac{t^{\gamma}}{2^{\gamma/2}\sqrt{\pi}} \int_0^1\frac{q^{\gamma/2 -1}}{\sqrt{1-q}}  i_{\gamma/2}(t\sqrt{1-2zq})dq = \frac{t^{\gamma}}{2^{\gamma/2}}K_{\gamma}(t,z).
\end{align*} 
With regard to \eqref{Int1}, \eqref{Int2111} and \eqref{Int31}, one gets the integrand of \eqref{ObF} for all strictly positive integers $\gamma$ which we shall extend in the following subsection to strictly positive real values of $\gamma$.

\subsection{Extension to $\Re(\gamma)> 0$} \eqref{ObF} extends to all complex values of $\gamma$ lying in the open right half-plane in a similar way as in remark p.194 in \cite{Dij}. The extension needs an exponential growth of both sides of \eqref{ObF} viewed as functions of the variable $\gamma$. To see this, we start with rewriting the LHS of \eqref{ObF} as 
 \begin{equation*}
\Gamma(\gamma)\int_{-1}^1 \sum_{s=1}^4\sum_{j \geq 0}(j+ \gamma)i_{j +\gamma}(\rho r) C_{j}^{\gamma}\left(c_{2\phi,2\theta}(u,v) w\right)
\left(\frac{\rho r}{2}\right)^j\xi^{js}\mu^{(\gamma-1)/2}(dw),
\end{equation*}
then, we use the bound  (\cite{Erd} p.14)
\begin{equation*}
|j_{\nu}(z)| \leq \frac{e^{|\Im(z)|}}{\Gamma(\nu+1)}, \quad z \in \mathbb{C}. 
\end{equation*}
which, together the definition of $i_{\nu}$, yields
\begin{equation*}
|i_{\nu}(z)| \leq \frac{e^{|z|}}{\Gamma(\nu+1)} \quad \Rightarrow \quad |\nu i_{\nu}(z)| \leq \frac{e^{|z|}}{\Gamma(\nu)} 
\end{equation*}
for $z \in \mathbb{R}$ and positive large enough $\nu$. Thus, the LHS of \eqref{ObF} is bounded by 
\begin{equation*}
4e^{\rho r}\int_{-1}^1\sum_{j \geq 0} |C_{j}^{\gamma}\left(c_{2\phi,2\theta}(u,v) w\right)| \frac{1}{(\gamma)_j}\left(\frac{\rho r}{2}\right)^j\mu^{(\gamma-1)/2}(dw).
\end{equation*}
Using the bound \eqref{BF}, Stirling's formula and since $\mu^{(\gamma-1)/2}$ has unit mass, then the LHS of \eqref{ObF} is easily seen to be bounded by 
\begin{equation*}
4e^{\rho r}\sum_{j \geq 0}\frac{(2\gamma)_j}{(\gamma)_j} \frac{(\rho r)^j}{2^j j!} =  4e^{\rho r} {}_1F_1(2\gamma, \gamma, (\rho r/2))
\end{equation*} 
which is of order $O(e^{\gamma})$. Coming to the RHS of \eqref{ObF}, note that $\Gamma((\gamma+1)/2)i_{(\gamma-1)/2}(z)$ is bounded by $e^{|z|}$ (as a function of $\gamma$ for fixed $z \in \mathbb{R}$), thereby so are 
\eqref{Int1} and \eqref{Int2111}. Finally, since $t^{\gamma} = e^{\gamma \log t} \in ]0,1]$ for all $t \in [0,1]$ and again $\Gamma(\nu+1)i_{\nu}(z)$ is bounded by $e^{|z|}$, then one only needs to focus on 
\begin{equation*}
\frac{\Gamma((\gamma+1)/2)\Gamma((\gamma/2)+j)}{[\Gamma((\gamma+1)/2+j)]^2}
\end{equation*} 
for all fixed $j \geq 0$, which is even bounded as a function of $\gamma$ by Stirling's formula. Thus, \eqref{ObF} holds for all $\gamma > 0$. Finally, \eqref{ObF} follows after integrating with respect to the symmetric Beta distributions $\mu^{l_0}, \mu^{l_1}$ displayed in \eqref{Dunkl} and the proof of \eqref{ObF} is finished.  
\begin{nota}
The function $Z \mapsto \cosh(t\sqrt{1-Z^2})$ appeared in  $($\cite{Boz}$)$ in relation to the famous Bessis-Moussa-Villani's conjecture, where the author used the Dirac-Pauli matrices:
\begin{equation*}
Q_1 := \left(\begin{array}{cc}
0 &1 \\
1 & 0 \end{array}\right), \qquad
Q_2 := \left(\begin{array}{cc}
0 &i \\
-i & 0 \end{array}\right),
\end{equation*}
then it is easy to see that 
\begin{equation*}
\cosh(t\sqrt{1-Z^2}) = \mathop{\rm tr}[e^{t(Q_1+ itZQ_2)}].
\end{equation*}
Moreover, in the Lie group scope, $Q_1$ and $Q_2$ are the generators of the $B_2$-Weyl group of the compact symplectic group $($\cite{Dunkl0}$)$.\\
\end{nota}

\subsection{Extension to $\gamma=0$}
It is obvious that \eqref{Int1} and \eqref{Int2111} have no singularity at $\gamma=0$. Hence, one only needs to deal with \eqref{Int3} where the only trouble comes from the constant term $(j=0)$, more precisely from $\Gamma(\gamma/2)$ in $A_{\gamma}$ (note that $\partial_zA_{\gamma}$ does not contain that factor since the series starts from $j=1$). But, it is known that $\sin(\gamma \pi/2)\Gamma(\gamma/2)$ has a removable singularity at $\gamma=0$ since the mirror formula holds 
\begin{equation*}
\Gamma(z)\Gamma(1-z) = \frac{\pi}{\sin \pi z}, \quad 0 < z < 1.
\end{equation*}
Moreover, since $(\sin \gamma \pi/2)/\gamma$ is well defined at $\gamma=0$, then the integrand in \eqref{Int31} may be extended to $\gamma=0$ and even tends to zero as $\gamma$ does. Since $\sin(\gamma \pi/2)$ remains in front of the integrals after the extension, then \eqref{Int31} tends to $0$ as $\gamma$ does. As a result
\begin{align*}
\lim_{k_1,k_0 \rightarrow 0^+}D_k^W(x,y) =  & \frac{1}{2}\lim_{k_1,k_0 \rightarrow 0^+} \Gamma\left(\frac{\gamma+1}{2}\right)
\int_{-1}^1\int_{-1}^1\left[ i_{(\gamma-1)/2}[\rho rc_{2\phi,2\theta}(u,v)] \right. \\&+ \left. \cos^2(\gamma\pi/2)i_{(\gamma-1)/2}[\rho r\sqrt{1-c_{2\phi,2\theta}^2(u,v)}]\right] 
\mu^{l_1}(du)\mu^{l_0}(dv).
\end{align*}
Since 
\begin{equation*}
1-c_{2\phi,2\theta}^2(u,v) = \frac{1-u\cos2\theta\cos2\phi - v \sin 2\theta \sin 2\phi}{2} = c_{2\phi,2\theta}^2(-u,-v),
\end{equation*}
then \eqref{Int1} and \eqref{Int21} give the same contribution up to the factor $\cos^2(\gamma\pi/2)$. Thus, 
\begin{align*} 
D_0^W(x,y) = \lim_{k_1,k_0 \rightarrow 0^+} \int_{-1}^1\int_{-1}^1 \Gamma\left(\frac{\gamma+1}{2}\right) i_{(\gamma-1)/2}[\rho rc_{2\phi,2\theta}(u,v)]\mu^{l_1}(du) \mu^{l_0}(dv) .
\end{align*}
(in fact, $i_{\nu}$ is an even function so that $i_{\nu}(|c_{2\phi,2\theta}(u,v)|) = i_{\nu}(c_{2\phi,2\theta}(u,v))$). Now, it is an easy exercice to see that the symmetric Beta distribution 
$\mu^{\nu}$ converges weakly as $\nu \rightarrow -1/2$ to the symmetric Bernoulli distribution (use the integral representation of the spherical Bessel function $j_{\nu}$)
\begin{equation*}
\eta(du) := \frac{1}{2}\left[\delta_{-1}(du) + \delta_1(du)\right].
\end{equation*}
It follows from the integral representation of $i_{(\gamma-1)/2}$ that
\begin{align*}
D_0^W(x,y) &=  \lim_{k_1,k_0 \rightarrow 0^+} \int\int \int \cosh(\rho rc_{2\phi,2\theta}(u,v)z)\mu^{(\gamma-1)/2}(dz)\mu^{l_1}(du) \mu^{l_0}(dv)
\\& =  \lim_{k_1,k_0 \rightarrow 0} \sum_{m=0}^{\infty}\frac{1}{(2m)!}\int\int \int \left[\frac{1+u\cos2\theta\cos2\phi + v\sin2\theta\sin 2\phi}{2}\right]^m(\rho rz)^{2m} \\& \mu^{(\gamma-1)/2}(dz)\mu^{l_1}(du) \mu^{l_0}(dv)
\\& = \int\int \int \cosh(\rho rc_{2\phi,2\theta}(u,v)z)\eta (dz)\eta(du) \eta(dv)
 \\& =  \int \int \cosh[\rho rc_{2\phi,2\theta}(u,v)]\eta(du) \eta(dv) 
\end{align*}
where we used Tonelli's Theorem to exchange the order of integration coming to the second equality, we used the fact that the weak convergence is equivalent to the convergence of moments for compactly-supported distributions together with the boundedness of the moments of $\mu^{\nu}$ by $1$ and Lebesgue's convergence Theorem to derive the third equality. Since 
\begin{equation*}
c_{2\phi,2\theta}(u,v) = \sqrt{\frac{1+z_{2\phi,2\theta}(u,v)}{2}},\, z_{2\phi,2\theta}(u,v) = u\cos 2\theta\cos2\phi + v \sin2\theta\sin2\phi,
\end{equation*}
and since $\cosh$ is an even function, then elementary trigonometric identities 
\begin{equation*}
\frac{1+\cos (2z)}{2} = \cos^2z,\quad \frac{1-\cos (2z)}{2} = \sin^2z
\end{equation*}
yield
\begin{align*}
&D_0^W(x,y)  = \frac{1}{4}\left[\cosh(\rho r \cos(\phi + \theta)) + \cosh(\rho r \cos(\phi - \theta)) + \cosh(\rho r \sin(\phi + \theta)) + \cosh(\rho r \sin(\phi - \theta))\right]
\end{align*}
which fits \eqref{GBF0} already derived in the introductory section.

\section{Proof of Proposition \eqref{Prop2}}
Let $\k,m$ be positive integers, $|\cdot|$ be the Euclidean norm in $\mathbb{R}^2$ and $Y_m$ is the ordinary spherical harmonic of degree $m$ written in polar coordinates as  
(see \cite{X98} for instance)
\begin{equation*}
Y_m (y) = |y|^m \cos m\t, \quad y= |y|e^{i\t}. 
\end{equation*}
Recall also that the generalized Bessel function is defined by (\cite{Dunkl1})
\begin{equation*}
D_k^W(x,y) := \frac{1}{8} \sum_{w \in B_2}D_k(x,wy) = \frac{1}{8} \sum_{w \in B_2} V_k \left[e^{\la \cdot, y w\ra} \right](y), 
\end{equation*}
where $D_k, V_k$ denote the Dunkl kernel and the Dunkl intertwining operator respectively (\cite{Dunkl1}). With the above notation, \eqref{D1} reduces with $p=2$ to
\begin{equation*}
D_0^W(x,y)= \frac{1}{8} \sum_{w \in B_2} e^{\la x, y w\ra}  =  I_0(\rho |y|) + 2  \sum_{j=1}^\infty \frac{ I_{4j}(\rho |y|)}{|y|^{4j} } Y_{4j}(y) \cos (4j \p).
\end{equation*}
Consequently, by linearity of $V_k$ (\cite{Dunkl1}), one has 
\begin{equation} 
\label{GBF1}
D_k^W(x,y) = V_k  \left[I_0(\rho |\cdot|)\right](y) + 2  \sum_{j=1}^\infty V_k \left[ \frac{ I_{4 j}(\rho |\cdot|)}{|\cdot|^{4j} } Y_{4j}\right](y) \cos (4j \p).
\end{equation}
Clearly we can expand $p_j^{l_1,l_0}(\cos \p)$ as a sum of $\cos m \phi, 0 \leq m \leq j$ so that, 
\begin{equation*}
p_j^{l_1,l_0}(\cos 4 \p) = \sum_{m=0}^j b_{m,j} \cos (4 m \p) 
\end{equation*}
(the coefficients $b_{m,j}$ may not be simple in general, although it has nice formula in the case of $l_1 = l_0$, that is, for Gegenbauer polynomials). 
Substituting this expression into \eqref{GBF} taken with $p=2$ and changing the summation's order, we get 
\begin{equation*}
D_k^W(x,y) =  c_{2,k}\left(\frac{2}{\rho r}\right)^{\gamma}\sum_{m= 0}^\infty \left( \sum_{j=m}^\infty b_{m,j}  I_{4j +\gamma}(\rho r)p_{j}^{l_1, l_0}(\cos(4\t)) \right)\cos(4m\p).
\end{equation*}
Comparing the last equation with (\ref{GBF1}), we obtain upon using the orthogonality of $\cos 4m\p$ that 
\begin{align}\label{Vn}
V_k \left[ \frac{ I_{4m}(\rho|\cdot|)} {(\rho|\cdot|)^{4m} } Y_{4m}\right](y) = \frac{2^{\gamma} c_{2,k}}{\rho^{4m}} \sum_{j=m}^\infty b_{m,j}  
\frac{ I_{4j +\gamma}(\rho|y|)} {(\rho |y|)^{\gamma}}p_{j}^{l_1, l_0}(\cos(4\t)).  
\end{align}

Using the expansion of the modified Bessel function \eqref{BessExp} to expand $I_{4j+\g}$ in \eqref{Vn}, the RHS of \eqref{Vn} becomes 
\begin{align*}
  &  \frac{2^{\gamma} c_{2,k}}{\rho^{4m}} \sum_{j=m}^\infty b_{m,j}  \sum_{q=0}^\infty\frac{1} {q! \Gamma(4 j + \g + q + 1)} \left(\frac{\rho|y|}{2}\right)^{2q+4j}
  p_{j}^{l_1, l_0}(\cos(4\t)) \\
  & =   c_{2,k}  \sum_{q = 0}^\infty \sum_{j=0}^\infty b_{m,j+m}   \frac{\rho^{4j+2q}}  {q! \Gamma(4 j +4m+ \g + q + 1)}
 \left(\frac{|y|}{2}\right)^{2q+4j + 4m} p_{j+m}^{l_1, l_0}(\cos(4\t)) \\      
  & =    c_{2,k}  \sum_{\k = 0 }^\infty \left[\sum_{q+2j=\k}  \frac{ b_{m,j+m}} {q! \Gamma(4 j +4m + q + \g +1)}
                 \left(\frac{|y|}{2}\right)^{2\k + 4m} p_{j+m}^{l_1, l_0}(\cos(4\t)) \right] \rho^{2\k}.            
\end{align*}
Similar computations transform the Bessel function in the LHS of \eqref{Vn} to
\begin{equation*}
\frac{ I_{4m}(\rho|\cdot|)} {(\rho|\cdot|)^{4m}}  = \frac{1}{2^{4m}}\sum_{\k=0}^{\infty}\frac{1} {\k! \Gamma(4m + \k +1)} \left(\frac{|\cdot|}{2}\right)^{2\k} \rho^{2\k}.
\end{equation*}
Comparing the coefficients of the obtained expressions, one finally gets 
\begin{equation*}
\frac{V_k \left[|\cdot |^{2\k} Y_{4m} \right](y)}{\k!\Gamma(4m+\k+1)}  = \sum_{0 \le 2j \le \k}  \frac{ b_{m,j+m}  |y|^{2\k+4m} } {(\k-2j)! \Gamma(4m + 2j+ \k+ \g +1)}p_{j+m}^{l_1, l_0}(\cos(4\t))  
\end{equation*}
and Proposition \eqref{Prop2} is proved.

{\bf Acknowledgments}: The author wants to thank deeply Professor Y. Xu from University of Oregon who informed him how to derive Proposition \eqref{Prop2} from \eqref{GBF}. 
He also wishes to thank Professor T. H. Koornwinder for providing him with important references . Special thanks to Professor C. F. Dunkl from UVa for his contribution to a better organization of the manuscript.

\end{document}